\documentclass[a4paper, 12pt]{amsart}
\usepackage{amsmath, amsthm, amscd, amssymb, amsfonts, amsxtra, amssymb, latexsym, bm}
\usepackage{enumerate}
\usepackage{verbatim}
\usepackage{float}
\usepackage{hyperref}
\hypersetup{colorlinks = true,	allcolors  = blue}

\newtheorem{thm}{Theorem}[section]
\newtheorem{prop}[thm]{Proposition}
\newtheorem{lem}[thm]{Lemma}
\newtheorem{coro}[thm]{Corollary}

\theoremstyle{definition}
\newtheorem{rem}[thm]{Remark}
\newtheorem{exam}[thm]{Example}

\theoremstyle{remark}

\usepackage{color}

\begin{document}
\numberwithin{equation}{section}
\title[Nilpotent Lie algebras of nilpotence class at most $3$ via graphs]{Nilpotent Lie algebras of nilpotence class at most three via graphs}
\author[Andrea L. Gallo, Denis E. Videla]{Andrea L. Gallo, Denis E. Videla}
\dedicatory{\today}
\thanks{Partially supported by CONICET and SECyT-UNC}

\address{Andrea L.\@ Gallo. FaMAF -- CIEM (CONICET), Universidad Nacional de C\'ordoba. 
	\newline	Av.\@ Medina Allende 2144, Ciudad Universitaria, (5000), C\'ordoba, Argentina.
	\newline {\it E-mail: andrea.gallo@unc.edu.ar}}

\address{Denis E.\@ Videla. FaMAF -- CIEM (CONICET), Universidad Nacional de C\'ordoba. 
\newline	Av.\@ Medina Allende 2144, Ciudad Universitaria, (5000), C\'ordoba, Argentina.
\newline {\it E-mail: devidela@unc.edu.ar}}

\begin{abstract}
In this work, we associate with each finite simple graph $G$ a nilpotent Lie algebra $\mathfrak{n}_G$ of class at most three, whose underlying
vector space has a basis indexed by the vertices, edges, and paths of length two of $G$..
Its nilpotency class is exactly three if and only if $G$ contains a path of length two.
We compute the dimensions of its lower central series and describe its center. Over fields of characteristic different
from two, we prove that two such Lie algebras are isomorphic if and only if their underlying graphs are isomorphic.
We also give an algebraic description using quotients involving exterior and symmetric squares. Finally, we study balanced
blow-ups and prove that $\mathfrak{n}_G$ is a retract of $\mathfrak{n}_{G^{(m)}}$.
 \end{abstract}

\maketitle

\section{Introduction}

 In 2005, Dani and Mainkar \cite{DM} gave a construction of a nilpotent Lie algebra of class at most two from the vertices and edges of a finite simple graph. 
The vertex space generates the algebra, and the bracket of two distinct vertices is the corresponding edge vector when they are adjacent, and zero otherwise. 
These Lie algebras and their associated Lie groups have been studied from several geometric and algebraic perspectives (see \cite{AV}, \cite{La}, \cite{LW}, \cite{Ov}, \cite{PT})  

Several constructions extend this correspondence. 
In \cite{Mainkar2006}, Mainkar considered $k$-step nilpotent Lie algebras obtained as quotients of free nilpotent Lie algebras by relations determined by the non-edges of a graph. 
In \cite{GGI}, Grantcharov, Grantcharov and Iliev constructed solvable Lie algebras by incorporating additional generators associated with $3$-cliques. 
In \cite{Ra}, Ray associated two-step nilpotent Lie algebras with Schreier graphs and characterized when these algebras admit three-step nilpotent extensions.

In this paper, we introduce a construction based on the vertices, edges, and paths of length two of an arbitrary finite simple graph $G$. 
Its underlying vector space is
$$
\mathfrak n_G=V\oplus A\oplus W,
$$
where the three summands have bases indexed by these respective combinatorial objects. 
The bracket satisfies
$$
[V,V]\subseteq A,\qquad [V,A]\subseteq W,
$$
with $W$ central and $[A,A]=0$. 
We prove that these relations define a nilpotent Lie algebra of class at most three. 
Its class is exactly three if and only if $G$ contains a path of length two.

The construction retains a direct connection with the
Dani--Mainkar algebra $\mathfrak n_G^{(2)}$, since
$$
\mathfrak n_G/W\cong\mathfrak n_G^{(2)}
$$
for every graph $G$. 
In particular, the two constructions coincide when $G$ has no paths of length two. 
The contribution of a triple of vertices depends on its induced edges: an
induced path contributes one dimension to both $W$ and
$[V,A]$, whereas a triangle contributes three dimensions to
$W$ and two to $[V,A]$. 
This distinction allows us to give explicit formulas for the dimensions of the lower central series, 
in terms of vertex degrees and the number of triangles.

Our main classification result states that 
$$
\mathfrak n_G\cong\mathfrak n_H
\quad\Leftrightarrow\quad
G\cong H,
$$
for arbitrary finite simple graphs $G$ and $H$ and any field of characteristic different from two.

We also describe the construction using quotients of
$\Lambda^2V$ and $V\otimes\operatorname{Sym}^2V$. 
This description gives a uniform formula for the bracket.
Finally, we study balanced blow-ups of graphs: if $G^{(m)}$ is
obtained by replacing each vertex of $G$ with $m$ independent
copies, we see that $\mathfrak n_G$ is a retract of $\mathfrak n_{G^{(m)}}$. 
We compute the dimensions of the resulting algebras and their lower central series.

The paper is organized as follows. Section~2 introduces the
construction and establishes the general dimension formulas.
Section~3 treats several families of graphs, including trees,
cycles, and strongly regular graphs. Section~4 studies the
center and disjoint unions and proves the isomorphism theorem.
Section~5 gives the tensor description, and Section~6 treats
balanced blow-ups.

\medskip
\noindent\textbf{Graph conventions.}
Throughout the paper, all Lie algebras are defined over a field
$\Bbbk$, and all graphs are finite, simple, and undirected.
We write $G=(S,E)$ and denote the edge $\{v,w\}$ by $vw$.
A path of length $t$ from $v$ to $w$ is a sequence of different vertices
$v_0,\ldots,v_t$ such that $v_iv_{i+1}\in E$ for every
$i=0,\ldots,t-1$ with $v_0=v$ and $v_t=w$. 
Paths of length two are identified with their reversals. 
In particular, the endpoints of such a path may be adjacent.

For $v\in S$, the set $ N(v)=\{w\in S:vw\in E\}$ is called the neighborhood of $v$ and $d(v)=|N(v)|$ the degree of $v$.
Recall that a graph is $k$-regular if every vertex has degree $k$. 
A connected non-complete graph is strongly regular with parameters $(n,k,e,d)$ if it has $n$ vertices, is $k$-regular,
and $|N(v)\cap N(w)|$ is equal to $e$ or $d$ depending on $v$ and $w$ are adjacent or not, respectively.
A tree is an acyclic connected graph without cycles. 
We refer to \cite{Di} for graph theory terminology and to \cite{BrW} for strongly regular graphs.

\section{New $3$-step nilpotent Lie algebras from graphs}

Let $G=(S,E)$ be a finite simple graph with
$S=\{v_1,\ldots,v_n\}$. Let $V$ be the $\Bbbk$-vector
space with basis $S$, and let $A$ have basis
$$
\{a_{ij}:i<j,\ v_i v_j\in E\}.
$$
Let $W$ have basis indexed by the paths of length two in $G$, identified with their reversals.
 We write $w_{ijk}$ for the basis vector corresponding to the path $v_i v_j v_k$, with $i<k$.

We define a bracket on $\mathfrak{n}_G=V\oplus A\oplus W$.
For $i<j$, set 
$$
[v_i,v_j]=
\begin{cases}
a_{ij} &\quad \text{if } v_i v_j\in E, \\
0 & \quad \text{if } v_i v_j\not\in E.
\end{cases}
$$
For $i<j<k$, the brackets involving one vertex and one edge are defined as follows.
\begin{enumerate}[$(a)$]
	\item If $v_i v_j, v_i v_k, v_j v_k \in E$, i.e.\@ the vertices $v_i,v_j,v_k$ determine a triangle in $G$ then
	\begin{itemize}
		\item $[v_i,a_{jk}]=w_{ijk}-w_{ikj}$.
		\item $[a_{ij},v_k]=w_{ijk}-w_{jik}$.
		\item $[a_{ik},v_j]=w_{ikj}-w_{jik}$.
	\end{itemize}
	
	\medskip
	
	\item If the graph induced by $v_{i}, v_j, v_k$ in $G$ is a $2$-path, then 
	\begin{itemize}
		\item If $v_j v_k \in E$, then 
				$$[v_i,a_{jk}]= 
				\begin{cases}
					w_{ijk} & \text{if }v_i v_j\in E, \\
					-w_{ikj} & \text{if }v_i v_k \in E.
				\end{cases}
				$$
		\item If $v_i v_j \in E$, then 
		$$[a_{ij},v_{k}]= 
		\begin{cases}
			w_{ijk} & \text{if }v_j v_k\in E, \\
			-w_{jik} & \text{if }v_i v_k \in E.
		\end{cases}
		$$
		\item If $v_i v_k \in E$, then 
		$$[a_{ik},v_{j}]= 
		\begin{cases}
			w_{ikj} & \text{if }v_j v_k\in E, \\
			-w_{jik} & \text{if }v_i v_j \in E.
		\end{cases}
		$$
		\end{itemize}
\end{enumerate}

All brackets not specified above are zero. 
We have the following result that states that the $\Bbbk$-vector space
$$
\mathfrak n_{G}=V\oplus A \oplus W
$$ 
with brackets as above is a nilpotent Lie algebra. 

\begin{thm}\label{Teo: Lie alg graph}
	Let $G=(S,E)$ be a simple undirected graph and let $\ell(G)$ be the number of different $2$-paths in $G$. 
	Then, the $\Bbbk$-vector space $\mathfrak n_G$ with brackets defined as above is a nilpotent Lie algebra of class at most $3$ of dimension $|S|+|E|+\ell(G)$.
\end{thm}
\begin{proof}
	We first check that the prescribed brackets define a bilinear
	antisymmetric operation on $\mathfrak n_G=V\oplus A\oplus W$. 
	We extend the prescribed brackets by bilinearity, impose $[b,b]=0$ for every basis element $b$, and define $[b',b]=-[b,b']$ for distinct basis elements. 
	All unspecified brackets are zero. The resulting operation is alternating,
	including in characteristic two. 
	In particular,
	$$
	[A,A]=0,\qquad [\mathfrak n_G,W]=0, \qquad [V,V]\subseteq A,\qquad [V,A]\subseteq W.
	$$
	
	We now verify the Jacobi identity. The preceding inclusions show
	that a Jacobi expression involving an element of $A$ or $W$
	vanishes: if one of its entries belongs to $W$, this follows from
	$[\mathfrak n_G,W]=0$; if one entry belongs to $A$ and the other two belong
	to $V$, its three terms vanish because $[A,A]=[V,W]=0$.
	The remaining cases with at least two entries in $A$ vanish for
	the same reason. By trilinearity, it is therefore enough to
	consider three vertices.
	
	If two vertices coincide, the Jacobi identity follows from alternation and the relation $[x,y]=-[y,x]$. 
	Let $i<j<k$ and consider the three distinct vertices $v_i,v_j,v_k$. Their Jacobi expression is
	$$
	J(v_i,v_j,v_k) =[v_i,[v_j,v_k]] +[v_j,[v_k,v_i]] +[v_k,[v_i,v_j]].
	$$
	We distinguish the possibilities for the subgraph induced by $v_i,v_j,v_k$.
	
	If it has at most one edge, each term in $J(v_i,v_j,v_k)$ is zero. 
	If it has exactly two edges, it is a $2$-path. There are three cases.
	
	If $v_iv_j,v_jv_k\in E$ and $v_iv_k\notin E$, the defining
	relations give
	$$
	[v_i,[v_j,v_k]]=w_{ijk},\qquad
	[v_j,[v_k,v_i]]=0,\qquad
	[v_k,[v_i,v_j]]=-w_{ijk}.
	$$
	If $v_iv_k,v_jv_k\in E$ and $v_iv_j\notin E$, they give
	$$
	[v_i,[v_j,v_k]]=-w_{ikj},\qquad
	[v_j,[v_k,v_i]]=w_{ikj},\qquad
	[v_k,[v_i,v_j]]=0.
	$$
	Finally, if $v_iv_j,v_iv_k\in E$ and $v_jv_k\notin E$, they give
	$$
	[v_i,[v_j,v_k]]=0,\qquad
	[v_j,[v_k,v_i]]=-w_{jik},\qquad
	[v_k,[v_i,v_j]]=w_{jik}.
	$$
	Thus, we have $J(v_i,v_j,v_k)=0$ in each of the three cases.
	
	If $v_i,v_j,v_k$ form a triangle, the defining relations yield
	\begin{align*}
		[v_i,[v_j,v_k]] &=w_{ijk}-w_{ikj},\\
		[v_j,[v_k,v_i]]&=w_{ikj}-w_{jik},\\
		[v_k,[v_i,v_j]]&=w_{jik}-w_{ijk}.
	\end{align*}
	These terms sum to zero, proving the Jacobi identity. Therefore $\mathfrak{n}_G$ is a Lie algebra.
	
	It remains to determine its nilpotency class. Every edge vector
	$a_{ij}$ is a bracket $[v_i,v_j]$, and all brackets of vertices
	with elements of $A$ lie in $W$. 
	Hence, the lower central series satisfies
	$$
	\gamma_2(\mathfrak n_G)=[\mathfrak n_G,\mathfrak n_G]=A\oplus[V,A],
	\qquad
	\gamma_3(\mathfrak n_G)=[\mathfrak n_G,\gamma_2(\mathfrak n_G)]=[V,A]\subseteq W,
	\qquad
	\gamma_4(\mathfrak n_G)=0.
	$$
	Thus $\mathfrak n_G$ is nilpotent of class at most three. 
	If $G$ contains a $2$-path, then at least one of the displayed brackets in the $2$-path case is non-zero,  so $[V,A]\ne0$ and the class is exactly three. 
	The same conclusion holds when $G$ contains a
	triangle, as its defining brackets are non-zero.
	
	Finally, $S$ is a basis of $V$, the edges form a basis of $A$,
	and the $2$-paths indexing the elements $w_{ijk}$ form a basis
	of $W$. Therefore
	$$
	\dim \mathfrak n_G=\dim V+\dim A+\dim W =|S|+|E|+\ell(G),
	$$
	as asserted.
\end{proof}	

%
%
%
%

\begin{rem}
	For a general graph $G$, the equality
	$[\mathfrak n_G,\mathfrak n_G]=A\oplus W$ need not hold.
	The bracket relations imply
	$$
	[\mathfrak n_G,\mathfrak n_G]
	=A\oplus[V,A]\subseteq A\oplus W,
	\qquad
	[\mathfrak n_G,[\mathfrak n_G,\mathfrak n_G]]
	=[V,A]\subseteq W.
	$$
	Therefore, the derived algebra is a subalgebra of the abelian
	subalgebra $A\oplus W$. The inclusions may be strict, as
	shown by the complete graph $K_3$ below.
\end{rem}

\begin{exam}
	Let $G=K_3$, with vertices $v_1,v_2,v_3$. Write
	$$
	A=\operatorname{span}_{\Bbbk}\{a_{12},a_{13},a_{23}\},
	\qquad
	W=\operatorname{span}_{\Bbbk}\{w_{123},w_{132},w_{213}\}.
	$$
	Here the three basis elements of $W$ correspond to the
	three paths of length two in the triangle, with the two
	endpoints ordered increasingly. Hence
	$$
	\mathfrak n_{K_3}=V\oplus A\oplus W,
	\qquad
	\dim\mathfrak n_{K_3}=3+3+3=9.
	$$
	The non-zero brackets between basis elements, up to
	antisymmetry, are
	\[
	[v_1,v_2]=a_{12},\qquad
	[v_1,v_3]=a_{13},\qquad
	[v_2,v_3]=a_{23},
	\]
	and
	$$	
	[v_1,a_{23}]=w_{123}-w_{132},\quad [a_{12},v_3]=w_{123}-w_{213},\quad [a_{13},v_2]=w_{132}-w_{213}.
	$$
	All remaining brackets between basis elements are zero.
	In particular, $W$ is central. 
	Notice that
	$$
	[V,A] =\operatorname{span}_{\Bbbk}\{ w_{123}-w_{132},\, w_{132}-w_{213}\},
	\qquad
	\dim[V,A]=2.
	$$
	Consequently, its lower central series is
	$$
	\gamma_2(\mathfrak n_{K_3}) =[\mathfrak n_{K_3},\mathfrak n_{K_3}] =A\oplus[V,A], \qquad
	\gamma_3(\mathfrak n_{K_3})=[V,A], \qquad
	\gamma_4(\mathfrak n_{K_3})=0.
	$$
	In particular, $\mathfrak n_{K_3}$ has nilpotency class
	exactly three, while
	$$
	\dim\gamma_2(\mathfrak n_{K_3})=5, \qquad
	\dim\gamma_3(\mathfrak n_{K_3})=2.
	$$
	Thus
	$$
	[\mathfrak n_{K_3},\mathfrak n_{K_3}]\subsetneq A\oplus W,
	$$
	so the three-dimensional space $W$ contains a one-dimensional central direction which is not generated by triple brackets. \hfill $\lozenge$
\end{exam}

\begin{prop}\label{prop: dimensiones}
	Let $G=(S,E)$ be a finite simple graph. Denote by $p(G)$
	the number of three-element subsets of $S$ inducing exactly
	two edges, and by $\tau(G)$ the number of triangles. 
	Then, $\dim\mathfrak n_G =|S|+|E|+p(G)+3\tau(G)$. 
	Furthermore, for the lower central series of $\mathfrak n_G$ we have
	$$
	\gamma_2(\mathfrak n_G)=A\oplus[V,A],
	\qquad
	\gamma_3(\mathfrak n_G)=[V,A],
	\qquad
	\gamma_4(\mathfrak n_G)=0,
	$$
	in this case
	$$
	\dim\gamma_2(\mathfrak n_G)=|E|+p(G)+2\tau(G), \qquad \dim\gamma_3(\mathfrak n_G)=p(G)+2\tau(G).
	$$
	Consequently, the dimensions of the successive quotients
	are
	\[
	\dim\bigl(\mathfrak n_G/\gamma_2(\mathfrak n_G)\bigr)
	=|S|+\tau(G),\qquad
	\dim\bigl(\gamma_2(\mathfrak n_G)/
	\gamma_3(\mathfrak n_G)\bigr)=|E|,
	\]
	and
	\[
	\dim\gamma_3(\mathfrak n_G)=p(G)+2\tau(G).
	\]
\end{prop}

\begin{proof}
	The relations defining $\mathfrak n_G$ give
	\[
	\gamma_2(\mathfrak n_G)=A\oplus[V,A],
	\qquad
	\gamma_3(\mathfrak n_G)=[V,A],
	\qquad
	\gamma_4(\mathfrak n_G)=0.
	\]
	It remains to compute the dimension of $[V,A]$.
	
	The basis of $W$ is partitioned according to the
	three-element vertex sets supporting its $2$-paths.
	The brackets associated with distinct vertex sets
	therefore have disjoint basis supports.
	
	If a vertex set induces exactly two edges, it contains
	one $2$-path. The defining relations show that its basis
	vector occurs as a bracket in $[V,A]$, so this vertex set
	contributes one dimension. If a vertex set
	$\{v_i,v_j,v_k\}$, with $i<j<k$, induces a triangle,
	its three path vectors are $w_{ijk},w_{ikj},w_{jik}$.
	The brackets involving these vertices span
	$$
	\operatorname{span}_{\Bbbk}\{ w_{ijk}-w_{ikj},\, w_{ijk}-w_{jik},\, w_{ikj}-w_{jik}\},
	$$
	which has dimension two. Vertex sets inducing at most
	one edge contribute nothing. Since the contributions
	from distinct vertex sets are independent, we obtain
	$$
	\dim[V,A]=p(G)+2\tau(G).
	$$
	
	Each induced $2$-path contributes one basis vector to $W$, while each triangle contributes three. 
	Thus, $\dim W=p(G)+3\tau(G)$. 
	Combining these identities with $\dim V=|S|$ and $\dim A=|E|$ proves all the stated formulas.
\end{proof}

\begin{coro}\label{coro: dimensiones}
	Let $G=(S,E)$ be a finite simple graph, let $d(v)$ denote
	the degree of $v\in S$, and let $\tau(G)$ be the number of
	triangles in $G$. Then
	$$
	\dim\gamma_2(\mathfrak n_G)=\tfrac{1}{2}\sum_{v\in S}\big(d(v)\big)^2-\tau(G)
	\quad	\text{and}
	\quad	\dim\gamma_3(\mathfrak n_G)=\sum_{v\in S}\binom{d(v)}{2}-\tau(G).
	$$
\end{coro}

\begin{proof}
	First, recall that 
	$$
	2|E|=\sum_{v\in S} d(v).
	$$
	On the other hand, each $2$-path has a unique middle vertex. Thus, if we choose its
	two incident edges at that vertex we obtain
	$$
	\sum_{v\in S}\binom{d(v)}{2}
	=p(G)+3\tau(G),
	$$
	since a triple inducing exactly two edges contributes
	one $2$-path, whereas a triangle contributes three.
	Substituting
	\[
	p(G)=\sum_{v\in S}\binom{d(v)}{2}-3\tau(G)
	\]
	into the formulas of the preceding proposition gives
	the result.
\end{proof}

\subsection*{Some remarks about the construction}

We conclude this section by comparing our construction with
two other constructions of Lie algebras from graphs.

\noindent{($i$)} Notice that our construction does not coincide with the solvable Lie algebra
construction of Grantcharov, Grantcharov, and Iliev
\cite{GGI}. This can already be seen for the complete graph
$K_3$. In our construction $\dim \mathfrak n_{K_3}=9$, moreover $\mathfrak n_{K_3}$ is nilpotent of class three.
On the other hand, in the construction of \cite{GGI}, one element is added for the
unique $3$-clique of $K_3$. The resulting Lie algebra therefore
has dimension $3+3+1=7$. The clique element acts non-trivially
and diagonally on the vertex and edge spaces, so this algebra
is not nilpotent. In particular, the two Lie algebras associated
with $K_3$ are not isomorphic.

\noindent{($ii$)} Our construction also differs from the higher-step graph Lie algebra construction of Mainkar \cite{Mainkar2006}.
For $G=K_3$, Mainkar's construction is the free three-step nilpotent Lie algebra on three generators, 
its homogeneous components have dimensions $3$, $3$, and $8$, respectively and so  its total dimension is $14$. 
In contrast, the algebra $\mathfrak n_{K_3}$ constructed here has dimension $9$, with $\dim W=3$ and $\dim[\mathfrak n_{K_3},[\mathfrak n_{K_3},\mathfrak n_{K_3}]]=2$.
Therefore, the two algebras associated with $K_3$ are not isomorphic.

\section{Lie algebra associated with known families of graphs}
 
We apply the general dimension formulas to several families of graphs, including trees, paths, stars, cycles, and strongly regular graphs.

\subsection*{Trees} Recall that a tree is a connected graph without cycles. 
In particular, any tree $T=(S,E)$ is free of triangles, 
and it satisfies $|S|=|E|+1$.

\begin{prop}\label{coro: dimensiones trees}
Let $T$ be a tree of order $n\geq3$. Then $\mathfrak{n}_T$ is a three-step nilpotent Lie algebra of
dimension $2n-1+\ell(T)$, where $\ell(T)$ denotes the number of paths of length two in $T$.  
Furthermore, in this case
$$
\dim\gamma_2(\mathfrak n_T)=\tfrac{1}{2}\sum_{v\in S}\big(d(v)\big)^2
\quad	\text{and}
\quad	\dim\gamma_3(\mathfrak n_T)=\sum_{v\in S}\binom{d(v)}{2}.
$$ 
\end{prop}
\begin{proof}
	It is a direct consequence of Proposition~\ref{prop: dimensiones} and Corollary~\ref{coro: dimensiones}, 
	by taking into account that $\tau(T)=0$ since trees are triangle-free.
\end{proof}

\begin{exam}
	There are exactly two non-isomorphic trees on four vertices, the path
	$P_4$ and the star $K_{1,3}$. Both have four vertices, three edges,
	and no triangles.
	
	Notice that the number of $2$-paths of graphs $P_4$ and $K_{1,3}$ are $2$ and $3$, respectively. 
	By taking into account that the degree sequences of $P_{4}$ and $K_{1,3}$ are $(2,2,1,1)$ and $(3,1,1,1)$ respectively, 
	we obtain that
	$$
	\dim\mathfrak n_{P_4}=9,\qquad \dim\gamma_2(\mathfrak n_{P_4})=5,\qquad \dim\gamma_3(\mathfrak n_{P_4})=2.
	$$
	and 
	$$
	\dim\mathfrak n_{K_{1,3}}=10,\qquad \dim\gamma_2(\mathfrak n_{K_{1,3}})=6,\qquad \dim\gamma_3(\mathfrak n_{K_{1,3}})=3.
	$$
	
	In both cases, $\dim V=4$, $\dim A=3$, and
	$\gamma_4(\mathfrak n_G)=0$. Moreover, since both trees are connected
	and have at least three vertices, $Z(\mathfrak n_G)=W$; hence their
	centers have dimensions $2$ and $3$, respectively. \hfill $\lozenge$
\end{exam}

\smallskip
\noindent\textit{Paths.}
Let $\mathfrak u_{n+1}$ be the Lie algebra of strictly upper
triangular $(n+1)\times(n+1)$ matrices, and let
$$
I_4=\operatorname{span}_{\Bbbk}\{E_{ij}:j-i\geq 4\}.
$$
This is an ideal of $\mathfrak u_{n+1}$. For $n\geq 3$, there is an
isomorphism
$$
\mathfrak n_{P_n}\cong\mathfrak u_{n+1}/I_4.
$$
Indeed, if the vertices of $P_n$ are $v_1,\ldots,v_n$ in path order,
the isomorphism is given on the natural basis by
$$
v_i\longmapsto E_{i,i+1}+I_4,\qquad a_{i,i+1}\longmapsto E_{i,i+2}+I_4,\qquad w_{i,i+1,i+2}\longmapsto E_{i,i+3}+I_4.
$$
The required bracket relations follow from
$$
[E_{i,i+1},E_{i+1,i+2}]=E_{i,i+2}, \qquad
[E_{i,i+1},E_{i+1,i+3}] =[E_{i,i+2},E_{i+2,i+3}] =E_{i,i+3}.
$$
In particular, $\mathfrak n_{P_3}\cong\mathfrak u_4$.

\smallskip
\noindent\textit{Stars.}
The star graph is the bipartite graph $G=K_{1,m}$ with $m\geq 2$. 
Denote its middle vertex by $c$, its leaves by $x_1,\ldots,x_m$, and put $a_i=[c,x_i]$. For $i<j$,
let $z_{ij}$ correspond to the path $x_i-c-x_j$. Besides
$[c,x_i]=a_i$, the non-zero brackets are
$$
[x_i,a_j]=[x_j,a_i]=z_{ij}
\qquad (i<j).
$$
There are $m+1$ vertices, $m$ edges, and $\binom{m}{2}$ paths of length two.
Hence, we have that
$$
\dim\mathfrak n_{K_{1,m}} =2m+1+\tbinom{m}{2},\quad 
\dim\gamma_2(\mathfrak n_{K_{1,m}})=m+\tbinom{m}{2},\quad
\dim\gamma_3(\mathfrak n_{K_{1,m}})=\tbinom{m}{2}.
$$
Its successive lower central quotients therefore have dimensions
$\bigl(m+1,m,\binom{m}{2}\bigr)$, and
$\dim Z(\mathfrak n_{K_{1,m}})=\binom{m}{2}$.

\subsection*{Strongly regular graphs}

The following result shows that the dimensions of $\gamma_{\ell}(\mathfrak n_G )$ with $\ell=1,2,3$ of a strongly regular graph can be given in terms of its srg-parameters.

\begin{prop}\label{coro: dimensiones srg}
	Let $G$ be a strongly regular graph with parameters $(n,k,e,d)$. Then $\mathfrak{n}_G$ has nilpotency class
	three, and
	$$
	\dim\mathfrak n_G=n+\frac{nk^2}{2},
	\qquad
	\dim\gamma_2(\mathfrak n_G)
	=\frac{nk(3k-e)}{6},
	\qquad
	\dim\gamma_3(\mathfrak n_G)
	=\frac{nk(3k-e-3)}{6}.
	$$	
\end{prop}
\begin{proof}
Since $G$ is $k$-regular, the number of its edges is
$|E|=nk/2$, and the number of $2$-paths is
$$
\sum_{v\in S}\binom{d(v)}{2}
=n\binom{k}{2}.
$$
Moreover, each edge lies in $e$ triangles, so
\[
\tau(G)=\frac{nke}{6}.
\]
The dimension formula in Theorem~\ref{Teo: Lie alg graph} therefore gives
$$
\dim\mathfrak n_G=n+\frac{nk}{2}+n\binom{k}{2}=n+\frac{nk^2}{2}.
$$
Substituting the same identities into Corollary~\ref{coro: dimensiones}
yields the formulas for $\gamma_2$ and $\gamma_3$.
Finally, a connected strongly regular graph of order
$n\geq 3$ contains a $2$-path, so its associated algebra has nilpotency class three.
\end{proof}

\begin{exam}
	Consider the strongly regular graphs on nine vertices.
	It is well-known that there are, up to isomorphism, two strongly regular graphs with nine vertices, the $(3\times3)$-grid $H(2,3)$ and the complete
	tripartite graph $K_{3,3,3}$. We compute the dimensions of their
	associated Lie algebras directly from Proposition~\ref{coro: dimensiones srg}.
	
	For the $(3\times3)$-grid $H(2,3)$, its vertex set is
	$
	S=\{(i,j):1\leq i,j\leq3\},
	$
	where distinct vertices are adjacent if they have the same
	first or the same second coordinate, in this case its parameters are $(n,k,e,d)=(9,4,1,2)$. 
	Hence, we obtain that
	$$
\dim\mathfrak n_G =81,\quad \dim\gamma_2(\mathfrak n_G)=66,\quad \dim\gamma_3(\mathfrak n_G)=48.
	$$
	Thus, the dimensions of the successive quotients of the lower central
	series are $(15,18,48)$.
	
	On the other hand, the graph $K_{3,3,3}$ has srg-parameters
	$(n,k,e,d)=(9,6,3,6)$. Applying Proposition~\ref{coro: dimensiones srg} again gives
	$$
	\dim\mathfrak n_{K_{3,3,3}}= 171,\quad \dim\gamma_2(\mathfrak n_{K_{3,3,3}})=135,\quad \dim\gamma_3(\mathfrak n_{K_{3,3,3}})=108.
	$$
	In this case, the dimensions of the successive quotients are $(36,27,108)$. \hfill $\lozenge$
\end{exam}

\subsection*{Cycles.}

Let $G=C_n$ with $n\geq4$, and index its vertices cyclically by $x_1,\ldots,x_n$. Write $a_i=[x_i,x_{i+1}]$, with indices taken modulo $n$. 
For each $i$, the path $x_i-x_{i+1}-x_{i+2}$ determines an element $w_i\in W$, and the non-zero brackets involving an edge can be written as
$$
[x_i,a_{i+1}]=[a_i,x_{i+2}]=w_i.
$$
Since every vertex has degree two and $C_n$ has no triangles,
Corollary~\ref{coro: dimensiones} implies
$$
\dim V=n,\qquad \dim A=n,\qquad \dim W =\sum_{i=1}^{n}\binom{2}{2}=n.
$$
Consequently,
$$
\dim\mathfrak n_{C_n}=3n,\qquad \dim\gamma_2(\mathfrak n_{C_n})=2n,\qquad \dim\gamma_3(\mathfrak n_{C_n})=n.
$$
Therefore, the successive lower central quotients have dimensions
$(n,n,n)$. Moreover, $Z(\mathfrak n_{C_n})=W$.


\section{Structural properties of $\mathfrak{n}_G$}

We first show that disjoint unions of graphs correspond to
direct sums of their associated Lie algebras.

\begin{lem}\label{lem: unions}
	Let $G_1$ and $G_2$ be finite simple graphs such that $V(G_1)\cap V(G_2)=\varnothing$, and let $G=G_1\sqcup G_2$ be their disjoint
	union. Then,
	$$
	\mathfrak n_G\cong
	\mathfrak n_{G_1}\oplus\mathfrak n_{G_2}.
	$$
\end{lem}

\begin{proof}
	For $i=1,2$, let $V_i$, $A_i$, and $W_i$ be the vector spaces associated with $G_i$ respectively. 
	Notice that every vertex of $G$ belongs to exactly one of $G_1,G_2$. 
	Moreover, an edge of $G$ lies entirely in exactly one of these graphs. 
	The same is true of a $2$-path: since its two edges share a vertex, and there are no edges between $G_1$ and $G_2$. 
	Hence, the bases used in the construction give
	$$
	V=V_1\oplus V_2,\qquad
	A=A_1\oplus A_2,\qquad
	W=W_1\oplus W_2.
	$$
	Consequently,
	$$
	\mathfrak n_G
	=(V_1\oplus A_1\oplus W_1)
	\oplus(V_2\oplus A_2\oplus W_2)
	$$
	as a vector space.
	
	The defining brackets within $V_i\oplus A_i\oplus W_i$ are precisely those of $\mathfrak n_{G_i}$ for $i=1,2$. 
	Notice that all brackets between the two summands vanish. 
	Indeed, vertices in different components are non-adjacent, and a bracket between a vertex and an edge vector can be non-zero only when the corresponding vertices form a $2$-path or a triangle.
	Neither configuration can use vertices from both components.
	Brackets involving $W$ are zero by definition.
	Therefore, the displayed vector-space decomposition is a decomposition into commuting Lie ideals, proving the assertion.
\end{proof}

\begin{prop}\label{prop: components}
	Let $G_1,\ldots,G_s$ be the connected components of a
	finite simple graph $G$. Then,
	$$
	\mathfrak n_G\cong
	\bigoplus_{i=1}^{s}\mathfrak n_{G_i}.
	$$
	such that $[g_i,g_j]=0$ for any $g_i \in \mathfrak n_{G_i}$ and $g_j\in \mathfrak n_{G_j}$ with $i\neq j$.
	In particular,
	\[
	Z(\mathfrak n_G)\cong
	\bigoplus_{i=1}^{s}Z(\mathfrak n_{G_i}).
	\]
\end{prop}
\begin{proof}
	Since $G=G_1\sqcup\cdots\sqcup G_s$, applying Lemma~\ref{lem: unions}  we obtain that
	$$
	\mathfrak n_G\cong
	\mathfrak n_{G_1}\oplus\cdots\oplus\mathfrak n_{G_s}.
	$$
	An element $(x_1,\ldots,x_s)$ of this direct sum commutes
	with every element if and only if each $x_i$ commutes with
	every element of $\mathfrak n_{G_i}$ for any $i=1,\ldots,s$. 
	Therefore, the center of the direct sum is the direct sum of the centers, as
	asserted.
\end{proof}

Notice that given a simple graph $G$, if $G_*$ is the disjoint union of the connected
components of $G$ having at least three vertices, then
\begin{equation}\label{eq: decom G}
G=G_*\sqcup rK_2\sqcup sK_1
\end{equation}

for uniquely determined non-negative integers $r$ and $s$. 
As a direct consequence we obtain the following result that give us explicitly the center of $\mathfrak{n}_G$.

\begin{coro}
Let $G$ be a simple graph and let $G^*$ be the disjoint union of its connected components with at least three vertices. 
If $G$ satisfies \eqref{eq: decom G} then
$$
Z(\mathfrak n_{G})\cong  Z(\mathfrak{n}_{G^*}) \oplus \Bbbk^{r+s}
$$
\end{coro}
\begin{proof}
	By Proposition~\ref{prop: components},
	$$
	\mathfrak n_G\cong
	\mathfrak n_{G_*}\oplus
	\mathfrak n_{K_2}^{\oplus r}\oplus
	\mathfrak n_{K_1}^{\oplus s}.
	$$
	The center of $\mathfrak n_{K_2}$ has dimension one, whereas $\mathfrak n_{K_1}$ is one-dimensional and abelian.
	Taking centers of the direct sum gives
	$$
	Z(\mathfrak n_G)\cong Z(\mathfrak n_{G_*})
	\oplus\Bbbk^{r+s},
	$$
	as desired.
\end{proof}

\subsection*{The isomorphism theorem}
We first identify the center of the algebra associated with a connected graph having at least three vertices.

\begin{lem}\label{lem:center-connected}
	Let $H$ be a connected finite simple graph with at least three vertices. 
	Then,
	\begin{equation}\label{eq:center-connected}
		Z(\mathfrak n_H)=W_H.
	\end{equation}
\end{lem}
\begin{proof}
By definition $W_H\subseteq Z(\mathfrak{n}_H)$. Conversely, let $x=v+a+w$ be central, with $v\in V_H$, $a\in A_H$, and $w\in W_H$. 
Since $H$ has no isolated vertices, the $A_H$-components of the brackets
$[x,v_j]$ show that $v=0$. Indeed, if $v\neq 0$ then 
$$
v=\sum_{v_i\in V(H)}c_{i} v_i 
$$
with $c_{\ell}\neq 0$ for some $\ell=1,\ldots,n$. Thus, since $v_\ell$ is non-isolated there exists $t\in \{1,\ldots,n\}$ such that $v_\ell v_t \in E(H)$. 
In this case, the coefficient of the edge vector corresponding to $v_\ell v_t$ in $[v,v_t]$ is $\pm c_\ell$.
Thus $[v,v_t]\neq 0$, a contradiction and therefore, $v$ must be equal to $0$.

On the other hand, since $a\in A_H=\langle a_e\rangle_{e\in E(H)}$ we obtain that 
$$
a=\sum_{e\in E(H)}c_ea_e, \qquad c_e\in \Bbbk. 
$$ 
Let us fix an edge $e=uz\in E(H)$.
Since $H$ is connected and has at least three vertices, 
there is a vertex $w\notin\{u,z\}$ adjacent to $u$ or $z$.
The vertices $w,u,z$ induce a $2$-path or a triangle, and the
defining relations give $[w,a_{uz}]\ne0$. 
Notice that in the bracket 
$$
[w,a]=\sum_{e\in E(H)}c_e[w,a_e]
$$ 
only $a_{uz}$ can contribute to the coordinates of $W_H$
corresponding to paths on the vertex set $\{w,u,z\}$:
the other possible edges on this set are incident to $w$,
and their brackets with $w$ are zero. 
Since $[w,a]=0$, we obtain $c_e=0$. This holds for every edge $e$, so
$a=0$ and \eqref{eq:center-connected} follows.
\end{proof}

We are now ready to prove our main theorem.

\begin{thm}\label{teo: isom}
	Assume that $\operatorname{char}\Bbbk\neq2$.
	Let $G_1$ and $G_2$ be finite simple graphs, and let
	$\mathfrak n_{G_1}$ and $\mathfrak n_{G_2}$ be the Lie algebras constructed above.
	Then,
	$$
	\mathfrak n_{G_1}\cong \mathfrak n_{G_2}
	\qquad\Leftrightarrow\qquad
	G_1\cong G_2.
	$$
\end{thm}
\begin{proof}
	A graph isomorphism induces a linear isomorphism on the vertex spaces and corresponding isomorphisms on the edge and $2$-path spaces, with the signs of the edge vectors determined by their orientations. 
	These maps preserve the defining brackets. 
	Hence, $G_1\cong G_2$ implies $\mathfrak n_{G_1}\cong \mathfrak n_{G_2}$.
	
	Conversely, let $G_*$ be the disjoint union of the connected
	components of $G$ having at least three vertices. Then
	\[
	G=G_*\sqcup rK_2\sqcup sK_1
	\]
	for uniquely determined non-negative integers $r$ and $s$.
	By Proposition~\ref{prop: components},
	\begin{equation}\label{eq:component-decomposition}
		\mathfrak n_G
		\cong
		\mathfrak n_{G_*}\oplus
		\mathfrak n_{K_2}^{\oplus r}\oplus
		\mathfrak n_{K_1}^{\oplus s}.
	\end{equation}
	
%
%
%
%
	By Lemma~\ref{lem:center-connected}, every connected component
	$H$ of $G_*$ satisfies $Z(\mathfrak n_H)=W_H$. Therefore,
	Proposition~\ref{prop: components} gives
	$$
	Z(\mathfrak n_{G_*})=W_{G_*}.
	$$

	Now, the Lie algebra $\mathfrak n_{K_2}$ is the three-dimensional
	Heisenberg algebra, with center $A_{K_2}$, whereas
	$\mathfrak n_{K_1}$ is one-dimensional and abelian. Taking
	the quotient of \eqref{eq:component-decomposition} by its
	center, and using \eqref{eq:center-connected}, we obtain
	\[
	\mathfrak n_G/Z(\mathfrak n_G)
	\cong
	\mathfrak n^{(2)}_{G_*}\oplus\Bbbk^{2r}
	\cong
	\mathfrak n^{(2)}_{Q},
	\qquad
	Q:=G_*\sqcup 2rK_1,
	\]
	where $\mathfrak n^{(2)}_Q$ is the two-step graph Lie algebra
	of Dani and Mainkar. Since $\operatorname{char}\Bbbk \neq2$, by Mainkar's isomorphism theorem
	\cite[Theorem~1.1]{M1}, the isomorphism class of
	$\mathfrak n_G/Z(\mathfrak n_G)$ determines $Q$ up to graph
	isomorphism.
	
	Every component of $G_*$ has at least three vertices, so
	the isolated vertices of $Q$ are precisely those in
	$2rK_1$. Thus $Q$ determines both $G_*$ and $r$.
	Finally, \eqref{eq:component-decomposition} gives
	\[
	s=\dim\mathfrak n_G-\dim\mathfrak n_{G_*}-3r.
	\]
	Therefore, the isomorphism class of $\mathfrak n_G$
	determines $G_*$, $r$, and $s$, and hence determines
	$G=G_*\sqcup rK_2\sqcup sK_1$ up to graph isomorphism.
\end{proof}

Notice that in the proof we obtain the following relation between the Lie algebra $\mathfrak{n}_G$ and the Dani-Mainkar's Lie algebra.

\begin{coro}
Let $G$ be a simple graph and let $G^*$ be the disjoint union of its connected components with at least three vertices such that
$$
G=G_*\sqcup rK_2\sqcup sK_1.
$$ 
Then, 
$$
\mathfrak{n}_G/Z(\mathfrak{n}_G)\simeq n^{(2)}_{Q} \qquad \text{and} \qquad \mathfrak{n}_G/W_G\simeq \mathfrak n_{G}^{(2)},
$$
where $\mathfrak n_{Q}^{(2)}$ is the 2-step nilpotent Lie algebra associated with the graph $Q=G_* \sqcup 2rK_1$.
\end{coro}

The following example shows that the dimensions of the lower central series and the center do not determine the isomorphism class of a graph Lie algebra.

\begin{exam}
	Consider the trees $T_1$ and $T_2$ on the vertex set
	$S=\{a,b,c,d,e,f\}$
	with edge sets
	$$
	E(T_1)=\{ac,bc,cd,de,ef\}, \qquad E(T_2)=\{ac,ab,cd,de,cf\}.
	$$
	Both trees have degree sequence $(3,2,2,1,1,1)$, and hence
	$$
	p(T_i)=\sum_{v\in S}\binom{d(v)}{2} =\binom{3}{2}+2\binom{2}{2}=5\qquad (i=1,2).
	$$
	Since each tree has six vertices and five edges, Proposition~\ref{prop: dimensiones}
	gives
	$$
	\dim\mathfrak n_{T_i}=6+5+5=16, \qquad \dim\gamma_2(\mathfrak n_{T_i})=5+5=10,\qquad \dim\gamma_3(\mathfrak n_{T_i})=5.
	$$
	Hence, both algebras have successive lower central quotients of
	dimensions $(6,5,5)$. Moreover, both trees are connected and so
	$$
	Z(\mathfrak n_{T_i})=W_{T_i}, \qquad \dim Z(\mathfrak n_{T_i})=5.
	$$
	
	Nevertheless, $T_1$ and $T_2$ are not isomorphic. Indeed, in each tree, $c$ is the unique vertex of degree three.
	It has two neighbors of degree one in $T_1$, but only one
	in $T_2$.
	Consequently, Theorem~\ref{teo: isom} implies
	$$
	\mathfrak n_{T_1}\not\cong\mathfrak n_{T_2}.
	$$
	Therefore, even within the family of trees, these numerical invariants do not classify the associated Lie algebras.
\end{exam}

\section{Algebraic tensor description of $\mathfrak{n}_G$}

  
In this section, we show that $\mathfrak{n}_G$ can be described in terms of certain exterior and symmetric algebras.

Let $G=(S,E)$ be a finite simple graph, and let $V$ be a $\Bbbk$-vector space with basis $S=\{v_1,\ldots,v_n\}$. Define
$$
U_G=\operatorname{span}_{\Bbbk}\{v_i\wedge v_j:\, v_i v_j\notin E\},\qquad A_G=\Lambda^2V/U_G.
$$
Write $x\odot y$ for the class of $x\otimes y$ in $\operatorname{Sym}^2V$. Let $R_G$ be the subspace of
$V\otimes\operatorname{Sym}^2V$ spanned by the basis tensors $v_j\otimes(v_i\odot v_k)$ satisfying at least one of the conditions
$$
i=k,\qquad  v_i v_j\notin E,\qquad v_j v_k\notin E,
$$
and put
$$
W_G=(V\otimes\operatorname{Sym}^2V)/R_G.
$$

\begin{thm}
	Let $G=(S,E)$ be a simple undirected graph and let $V,\,A_G$ and $W_G$ be as defined above. 
	Then, $\mathfrak n_G$ is isomorphic to
	$$
	V\oplus A_G\oplus W_G
	$$
	with the bilinear alternating bracket determined by
	\[
	[x,y]=x\wedge y+U_G,
	\qquad  
	[x,y\wedge z+U_G]
	=
	y\otimes(x\odot z)-z\otimes(x\odot y)+R_G,
	\quad x,y,z\in V,
	\]
	and
	\[
	[A_G,A_G]=0,\qquad
	[V\oplus A_G\oplus W_G,W_G]=0.
	\]
	Under this isomorphism, the natural edge and $2$-path vectors
	correspond to
	$$
	a_{ij}\longmapsto v_i\wedge v_j+U_G,
	\qquad
	w_{ijk}\longmapsto v_j\otimes(v_i\odot v_k)+R_G.
	$$
\end{thm}

\begin{proof}
	The surviving basis elements of $A_G$ correspond exactly to the edges of $G$. 
	Similarly, a basis tensor $v_j\otimes(v_i\odot v_k)$ survives in $W_G$ precisely when $i\neq k$ and both $v_i v_j$ and $v_j v_k$ are edges.
	Since the symmetric product identifies the two orders of the end-points, these surviving tensors correspond exactly to the $2$-paths indexing the basis of $W$.
	
	We first verify that the prescribed bracket between $V$ and $A_G$ is well defined. 
	Notice that the expression
	$$
	y\otimes(x\odot z)-z\otimes(x\odot y)
	$$
	is bilinear and alternating in $y,z$, so it defines a linear map on $\Lambda^2V$ for each fixed $x$. 
	If $v_j v_k\notin E$, then for every basis vertex $v_i$ both
	$$
	v_j\otimes(v_i\odot v_k) \quad\text{and}\quad v_k\otimes(v_i\odot v_j)
	$$
	belong to $R_G$. Thus, this map vanishes on $U_G$ and descends to $A_G$.
	
	For $x,y,z\in V$, the Jacobi expression is the class modulo $R_G$ of
	\[
	y\otimes(x\odot z)-z\otimes(x\odot y)+z\otimes(y\odot x)-x\otimes(y\odot z)+x\otimes(z\odot y)-y\otimes(z\odot x),
	\]
	which is zero because the symmetric product is commutative.
	Every Jacobi expression involving an element of $A_G$ or $W_G$
	also vanishes, since $[V,A_G]\subseteq W_G$, $W_G$ is central,
	and $[A_G,A_G]=0$. Hence, the prescribed operation is a Lie bracket.
	
	Finally, for $i<j<k$, the formulas yield
	$$
	[v_i,a_{jk}]=w_{ijk}-w_{ikj}, \qquad	[a_{ij},v_k]=w_{ijk}-w_{jik}, \qquad [a_{ik},v_j]=w_{ikj}-w_{jik},
	$$
	where each term that does not correspond to a $2$-path of $G$ is interpreted as zero. 
	These are exactly the defining relations of $\mathfrak n_G$, both for triangles and for induced $2$-paths.
	Brackets involving a repeated vertex vanish by the definition of $R_G$. 
	Therefore, the stated basis identifications define a Lie algebra isomorphism.
\end{proof}

We finish this section with the following remark about the above algebraic description.

\begin{rem}
	
\noindent $(i)$	The tensor description above is valid over an arbitrary field, including fields of characteristic two. 
	Here, the exterior and symmetric squares are understood as quotients of the tensor square, and the Lie bracket is required to be alternating. 
	Neither the construction nor the verification of the Jacobi identity requires division by any scalar.

\noindent $(ii)$	In this description, the first tensor factor records the middle vertex of a $2$-path, whereas the symmetric product records its unordered end-points.
 	Consequently,
	$$
	\dim W_G=\sum_{v\in S}\binom{\deg(v)}{2}.
	$$
	In particular, a triangle contributes three independent vectors to $W_G$, one for each choice of middle vertex.
\end{rem}

\section{Blow-up operation and Lie algebras}

Recall that given a simple graph $G$, the $m$-balanced blow-up $G^{(m)}$ of $G$ is the graph obtained by replacing each vertex $v$ by an independent set
$\{v^{(1)},\ldots,v^{(m)}\}$, and $v^{(\alpha)}$ is adjacent to
$w^{(\beta)}$ if and only if $vw\in E(G)$ (See for instance \cite{HHN} and \cite{VJP}). 
As shown by Andrada and Vera in \cite{AV}, these constructions are useful for finding \textit{adapted complex structures} on $\mathfrak{n}_{G}^{(2)}$.
Here we study the relation between $\mathfrak{n}_G$ and $\mathfrak{n}_{G^{(m)}}$. 
We first express the dimensions of the lower central series of $\mathfrak{n}_{G^{(m)}}$ in terms of combinatorial invariants of $G$. 

\begin{prop}\label{prop: blow-up dim}
	Let $G=(S,E)$ be a finite simple graph, let $N=|S|$ and $M=|E|$ and let $m\geq2$. Then,
	$$
	\dim\mathfrak n_{G^{(m)}}=mN+m^3\bigl(M+\ell(G)\bigr),\quad \dim\gamma_2(\mathfrak n_{G^{(m)}}) =m^3\bigl(M+\ell(G)-\tau(G)\bigr),
	$$
	$$	
	\dim\gamma_3(\mathfrak n_{G^{(m)}})=m^3\bigl(\ell(G)-\tau(G)\bigr)+m^2(m-1)M.
	$$
\end{prop}

\begin{proof}
	Let $\ell(G)$ be the number of $2$-path in $G$, recall that
	$$
	\ell(G)=\sum_{v\in S}\binom{d(v)}{2},
	$$
	and let $\tau(G)$ be the number of triangles of $G$. Each vertex
	$v$ gives rise to $m$ vertices of degree $m d(v)$ in $G^{(m)}$.
	Moreover, each edge has $m^2$ lifts and each triangle has $m^3$
	lifts, this implies that
	$$
	|V(G^{(m)})|=mN, \qquad |E(G^{(m)})|=m^2M,\qquad \tau(G^{(m)})=m^3\tau(G).
	$$
	To compute	$\ell(G^{(m)})$, notice that every $2$-path is uniquely determined by its middle vertex
	and an unordered pair of distinct neighbors of that vertex.
	For each $v\in S$, the blow-up contains $m$ copies of $v$,
	each having $m d(v)$ neighbors. Thus, each copy is the middle
	vertex of exactly $\binom{md(v)}{2}$ paths, and consequently
	$$
	\ell(G^{(m)})
	=m\sum_{v\in S}\binom{md(v)}{2}.
	$$
	Now, notice that
	$$
	\binom{md(v)}{2}=\frac{md(v)(md(v)-1)}{2}=\frac{m^{2}d(v)^2-md(v)}{2}= m^2\binom{d(v)}{2}+\binom{m}{2}d(v).
	$$
	Therefore, by taking into account that $\sum_{v\in S}d(v)=2M$, we obtain
	$$
	\begin{aligned}
		\ell(G^{(m)}) &=m^3\sum_{v\in S}\binom{d(v)}{2}+m\binom{m}{2}\sum_{v\in S}d(v)\\
		&=m^3\ell(G)+m^2(m-1)M.
	\end{aligned}
	$$
	By Proposition~\ref{prop: dimensiones}, $\dim\mathfrak n_G=N+M+\ell(G)$, and
	$$
	\dim\gamma_2(\mathfrak n_G)=M+\ell(G)-\tau(G), \qquad \dim\gamma_3(\mathfrak n_G)=\ell(G)-\tau(G).
	$$
	Applying these formulas to $G^{(m)}$ gives
	$$
	\begin{aligned}
		\dim\mathfrak n_{G^{(m)}}&=mN+m^2M+m^3\ell(G)+m^2(m-1)M\\
		&=mN+m^3\bigl(M+\ell(G)\bigr),
	\end{aligned}
	$$
	and
	$$
	\begin{aligned}
		\dim\gamma_2(\mathfrak n_{G^{(m)}})
		&=m^2M+m^3\ell(G)+m^2(m-1)M-m^3\tau(G)\\
		&=m^3\bigl(M+\ell(G)-\tau(G)\bigr).
	\end{aligned}
	$$
	Finally,
	$$
	\dim\gamma_3(\mathfrak n_{G^{(m)}})
	=m^3\bigl(\ell(G)-\tau(G)\bigr)+m^2(m-1)M.
	$$
	These are the claimed identities.
\end{proof}

The following theorem states that $\mathfrak{n}_G$ is a retract of $\mathfrak{n}_{G^{(m)}}$.

\begin{thm}\label{teo: blow-up}
	Let $G$ be a finite simple graph and let $m\geq2$.
	Then $\mathfrak{n}_G$ is a retract of $\mathfrak{n}_{G^{(m)}}$. 
	More precisely, there are Lie algebra homomorphisms
	$$
	\mathfrak n_G	\xrightarrow{\ s\ } \mathfrak n_{G^{(m)}}\xrightarrow{\ \pi\ }\mathfrak n_G
	$$
	such that $\pi\circ s=\operatorname{id}_{\mathfrak n_G}$. In particular,
	$$
	\mathfrak n_{G^{(m)}}\cong \ker(\pi)\rtimes\mathfrak n_G.
	$$
\end{thm}

\begin{proof}
	We first define $\pi$ on the natural basis of
	$\mathfrak n_{G^{(m)}}$. On vertices and edges, let
	$$
	\pi(v_i^{(\alpha)})=v_i,\qquad \pi\bigl(a_{(i,\alpha)(j,\beta)}\bigr)=a_{ij},
	$$
	where the orientations of the edge vectors are chosen consistently.
	For a basis vector corresponding to the $2$-path
	$v_i^{(\alpha)}v_j^{(\beta)}v_k^{(\gamma)}$, define
	$$
	\pi\bigl(w_{(i,\alpha)(j,\beta)(k,\gamma)}\bigr)
	=
	\begin{cases}
		w_{ijk},& i\neq k,\\
		0,& i=k.
	\end{cases}
	$$
	As usual, the two endpoints in the notation for $w_{ijk}$ are ordered
	according to the convention used to define the basis of $W$.
	
	We check that $\pi$ preserves brackets. The bracket of two vertices
	corresponding to different fibers maps to the bracket of their images,
	because adjacency between fibers is determined by adjacency in $G$.
	Two vertices in the same fiber commute, and their images also commute.
	
	It remains to consider a bracket between a vertex and an edge vector.
	If the three vertices involved belong to distinct fibers, their
	induced edges are exactly those of the corresponding three vertices
	of $G$. Thus the defining bracket relations, including the signs in
	the triangle case, are carried to the defining relations of
	$\mathfrak n_G$. If two vertices belong to the same fiber, the image
	of every resulting $2$-path vector is zero. The bracket of the
	corresponding vertex and edge in $\mathfrak n_G$ is also zero, since
	it involves the same vertex twice. All remaining brackets vanish in
	both algebras. Hence $\pi$ is a Lie algebra homomorphism.
	
	Now choose one copy of each vertex, say $v_i^{(1)}$. These vertices
	induce a subgraph of $G^{(m)}$ isomorphic to $G$. The identification
	$v_i\mapsto v_i^{(1)}$, together with the corresponding identifications
	of edge and $2$-path vectors, therefore defines an injective Lie
	algebra homomorphism
	\[
	s:\mathfrak n_G\longrightarrow\mathfrak n_{G^{(m)}}.
	\]
	By construction, $\pi\circ s$ is the identity on every basis vector
	of $\mathfrak n_G$. This proves the retraction. Since $\ker(\pi)$ is
	an ideal and $s(\mathfrak n_G)$ is a subalgebra complementary to it,
	the stated semidirect product decomposition follows.
\end{proof}

The following example shows that blow-up can increase the nilpotency class of the associated Lie algebra. 
\begin{exam}
	The Lie algebra $\mathfrak n_{K_2}$ is the three-dimensional Heisenberg algebra, with
	$$
	\dim\gamma_2(\mathfrak n_{K_2})=1, \qquad \gamma_3(\mathfrak n_{K_2})=0.
	$$
	Since $K_2^{(m)}=K_{m,m}$, Proposition \ref{prop: blow-up dim} implies
	$$
	\dim\mathfrak n_{K_{m,m}}=2m+m^3,\qquad \dim\gamma_2(\mathfrak n_{K_{m,m}})=m^3, \qquad \dim\gamma_3(\mathfrak n_{K_{m,m}})=m^2(m-1).
	$$
	Thus, for every $m\geq2$, the algebra $\mathfrak n_{K_{m,m}}$ has nilpotency class exactly three,
	while $\mathfrak n_{K_2}$ has class two. 
	Its successive lower central quotients have dimensions
	$$
	\bigl(2m,m^2,m^2(m-1)\bigr).
	$$
	Furthermore, $K_{m,m}$ is connected and triangle-free, so
	$$
	Z(\mathfrak n_{K_{m,m}})=W_{K_{m,m}} =\gamma_3(\mathfrak n_{K_{m,m}}).
	$$
	By Theorem~\ref{teo: blow-up}, each of these algebras admits
	$\mathfrak n_{K_2}$ as a retract. \hfill $\lozenge$
\end{exam}

%

\end{document}